\documentclass[centertags,12pt,a4paper]{amsart}

\usepackage{latexsym}
\usepackage{amssymb,amsmath}
\usepackage{amscd}

\makeatletter
\@addtoreset{equation}{section}
\makeatother

\newtheorem{theorem}{Theorem}[section]
\newtheorem{proposition}[theorem]{Proposition}
\newtheorem{definition}[theorem]{Definition}
\newtheorem{lemma}[theorem]{Lemma}

\newtheorem*{claim*}{Claim}
\theoremstyle{definition}
\newtheorem{example}[theorem]{Example}
\newtheorem*{question*}{Question}

\def\cB{\mathcal{B}}
\def\cA{\mathcal{A}}

\def\cT{\mathcal{T}}

\def\cF{\mathcal{F}}

\def\cK{\mathcal{K}}
\def\cI{\mathcal{I}}

\def\cK{\mathcal{K}}

\def\fqss{\mathbb{F}_{q'}}
\def\fqssqq{\mathbb{F}_{\sqrt q}}

\def\fq{\mathbb{F}_q}

\author[Giulietti]{M. Giulietti}
\author[Montanucci]{E. Montanucci}
\thanks{2000 {\em Math. Subj. Class.}: Primary 51E21, Secondary 51A30, 05B25, 05C70}
\thanks{This research was performed within the activity of GNSAGA of the
Italian INDAM, with the financial support of the Italian Ministry MIUR,
project ``Strutture geometriche, combinatorica e loro applicazioni'',
PRIN 2004-2005}

\title[Hyperfocused arcs]{On hyperfocused arcs in $PG(2,q)$}

\address{Dipartimento di Matematica Universit\`a di Perugia, 06123
Perugia, Italy}\email{giuliet@dipmat.unipg.it,
montanuc@dipmat.unipg.it}

\begin{document}

   \begin{abstract} A $k$-arc in a Dearguesian projective plane whose
secants meet some external line in $k-1$ points is said to be
hyperfocused. Hyperfocused arcs are investigated in connection
with a secret sharing scheme based on geometry due to Simmons. In
this paper it is shown that point orbits under suitable groups of
elations are hyperfocused arcs with the significant property of
being contained neither in a hyperoval, nor in a proper subplane.
Also, the concept of generalized hyperfocused arc, i.e. an arc
whose secants admit a blocking set of minimum size, is introduced:
a construction method is provided, together with the
classification for size up to $10$.
\end{abstract}

\maketitle

     \section{Introduction}\label{s1}

Hyperfocused arcs were introduced in connection with a secret
sharing scheme based on geometry due to Simmons \cite{SIM}. The
implementation of this scheme needs an arc in a Desarguesian
projective plane with the property that its secant lines intersect
some external line in a minimal number of points. Simmons only
considered planes of odd order, where  this minimal number equals
the number of points of the arc \cite{BK}. He introduced the term
sharply focused set for arcs satisfying the aforementioned
property. Sharply focused sets in Desarguesian projective planes
of odd order were classified by Beutelspacher and Wettl \cite{BW},
whose result was based on a previous paper by Wettl \cite{WE}.

In 1997 Holder \cite{HOL} extended Simmons's investigation to
Desarguesian planes of even order. In such planes  the secants of
an arc of size $k$ may meet an external line in only $k-1$ points,
yet the classification of arcs having this property seems to be an
involved problem. Holder used the term super sharply focused sets
for such arcs and gave some constructions for them.

In a recent paper \cite{CH}, Cherowitzo and Holder proposed the
term hyperfocused arc instead of super sharply focused set. They
provided the classification of small hyperfocused arcs, and
constructed new examples, one of which gave a negative answer to a
question raised by Drake and Keating \cite{DK} on the possible
sizes of a hyperfocused arc.

Some open problems were pointed out by Cherowitzo and Holder,
including the existence of hyperfocused arcs which are neither
contained in a proper subplane nor in a hyperoval. In this paper a
positive answer to this question is given. The main tool is the
investigation of the so-called translation arcs, i.e. arcs which
are point orbits under a group of elations. In Section \ref{s2} it
is shown that such arcs are hyperfocused, and it is proved that
sometimes they are contained neither in a hyperoval nor in a
proper subplane, see Theorem \ref{sopra}.

The concept of hyperfocused arc can be naturally extended to that
of generalized hyperfocused arc, that is an arc of size $k$ for
which there exists an external point set of size $k-1$ meeting
each of its secants. Recently, Aguglia, Korchm\'aros and Siciliano
\cite{AKS} proved that in Desarguesian planes of even order any
generalized hyperfocused arc is hyperfocused, provided that it is
contained in a conic. In Section \ref{s3} we provide a
construction of generalized hyperfocused arcs which are not
hyperfocused. Also, a classification of small generalized
hyperfocused arcs is proved using the graph-theoretic concept of
$1$-factorizations of a complete graph, see Section \ref{s4}.

\section{Definitions and Notation}

Let $PG(2,q)$ be the Desarguesian plane over $\fq$, the finite
field with $q$ elements. A $k$-arc $\cK$ in $PG(2,q)$ is a set of
$k$ points no three of which are collinear. Any line containing
two points of $\cK$ is said to be a {\em secant} of $\cK$. A {\em
blocking set of the secants of $\cK$} is a point set $\cB\subset
PG(2,q)\setminus \cK$ having non-empty intersection with each
secant of $\cK$. As the number of secants of $\cK$ is $k(k-1)/2$,
the  size of $\cB$ is at least $k-1$. If this lower bound is
attained, $\cB$ is said to be {\em of minimum size}. Also, $\cB$
is {\em linear} if it is contained in a line.

Arcs in $PG(2,q)$ admitting a linear blocking set of minimum size
of their secants are called {\em hyperfocused arcs}. As mentioned
in the Introduction, hyperfocused arcs exist only in $PG(2,q)$ for
$q$ even. Therefore in the whole paper we assume $q=2^r$.

Throughout, we fix the following notation. Let $(X_1,X_2,X_3)$ be
homogeneous coordinates for points in $PG(2,q)$, and let
$\ell_\infty$ be the line of equation $X_3=0$. Given a pair
$A=(a,b)$ in $\fq\times\fq$, denote ${\overline A}$ the point in
$PG(2,q)$  with coordinates $(a,b,1)$, and  ${\overline A}_\infty$
the point $(a,b,0)$. Also, let $\varphi_A$ be the projectivity
$$
\varphi_A:(X_1,X_2,X_3)\mapsto(X_1+a_1X_3,X_2+a_2X_3,X_3)\,.
$$
Clearly, $\varphi_A$ is an elation with axis $\ell_\infty$, and
conversely for any non-trivial elation $\varphi$ with axis
$\ell_\infty$ there exists $A\in \fq \times \fq$, $A\neq (0,0)$,
such that $\varphi=\varphi_A$.

Given an additive subgroup $G$ of $\fq\times \fq$, let $\cK_G(P)$
be the orbit of the point $P\in PG(2,q)\setminus \ell_\infty$
under the action of the group
$$
T_G:=\{\varphi_A\mid A \in G\}\,.
$$
Clearly, any two orbits $\cK_G(P)$ and $\cK_G(Q)$ with $P$, $Q$
$\in PG(2,q)\setminus \ell_\infty$ are projectively equivalent.
For brevity, write $\cK_G$ for $\cK_G(O)$, where $O=(0,0,1)$. Note
that
$$
\cK_G:=\{{\overline A}\mid A\in G\}\,.
$$
A $k$-arc in $PG(2,q)$ coinciding with $\cK_G(P)$ for some
additive subgroup $G\subset \fq\times \fq$ and some $P\in
PG(2,q)\setminus \ell_\infty$ will be called a {\em translation
arc}.

\section{Translation arcs}\label{s2}
The following proposition shows that any translation arc is a
hyperfocused arc.

\begin{proposition}\label{prima}
Let $\cK$ be a translation arc. Then there exists a blocking set
 of the secants of $\cK$ of minimum size which is contained in
$\ell_\infty$.
\end{proposition}
\begin{proof}
Let $G$ be an additive subgroup of $\fq\times \fq$ such that $\cK$
is projectively equivalent to $\cK_G$. To prove the assertion, it
is enough to show that every secant of $\cK$ meets $\ell_\infty$
in a point ${\overline C}_\infty$ for some $C\in G\setminus
\{(0,0)\}$. For $A,B\in G$, $A\neq B$, let $l_{AB}$ be the secant
of $\cK$ passing through ${\overline A}$ and ${\overline B}$. The
intersection point of $l_{AB}$ and $ \ell_\infty$ is
${\overline{(A+B)}}_\infty$. Then the claim is proved, as $A+B$ is
a non-zero element of $G$.
\end{proof}

According to Proposition \ref{prima} groups $G$ in both Examples
\ref{n1} and \ref{n2} provide examples of translation arcs
$\cK_G$.

\begin{example}\label{n1} (see \cite{DK}) For any additive subgroup $H$ of
$\fq$, let $G=\{(\alpha,\alpha^2)\mid \alpha \in H\}$.
\end{example}

\begin{example}\label{n2}
For $H$ any additive subgroup of $\fq$ and $i$ any positive
integer with $(i,r)=1$, let $G=\{(\alpha,\alpha^{2^i})\mid \alpha
\in H\}$. Note that the arc $\cK_G$ is contained in a translation
hyperoval (see \cite[Ch. 8]{HIR}).
\end{example}




The following result shows that any translation $k$-arc is either
complete in $PG(2,q)\setminus \ell_\infty$ (i.e. it is not
contained in any $(k+1)$-arc $\cK'\subset PG(2,q)\setminus
\ell_\infty$), or it is contained in a translation $2k$-arc.
\begin{proposition}\label{laterale}
Let $\cK_G$ be a translation $k$-arc in $PG(2,q)$. Assume that
there exists a point $\overline{A}\in PG(2,q)$ belonging to no
secant of $\cK_G$. Then the set $\cK':=\cK_{G}\cup
\varphi_A(\cK_G)$
 is a translation $2k$-arc.
\end{proposition}
\begin{proof}
Assume that $\overline{A_1}$, $\overline{A_2}$ and
$\overline{A_3}$ are three collinear points in $\cK'$. Clearly
neither $\cK_G$ nor $\varphi_A(\cK_G)$ can contain all of such
points. Also, as $\varphi_A$ is an involution we may assume
$\overline{A_1},\overline{A_2}\in \cK_G$,
$\overline{A_3}\in\varphi_A(\cK_G)$. Note that the elation
$\varphi:=\varphi_{A+A_3}$ acts on both $\cK_G$ and
$\varphi_A(\cK_G)$. Then, as
$\varphi(\overline{A_3})=\overline{A}$, the secant of $\cK_G$
through $\varphi(\overline{A_1})$ and $\varphi(\overline{A_2})$
contains $\overline{A}$, which is a contradiction. Hence $\cK'$ is
a $2k$-arc. It is actually a translation arc because
$\cK'=\cK_{G'}$, where $G'=G\cup (G+A)$.

\end{proof}


\begin{example}\label{n3}
For $q$ a square, let $\eta\in \fq\setminus \fqssqq$ such that
$\eta^2\in \fqssqq$. Choose $b\in \fqssqq$, $b\neq 1$, and let
$G=\{(\alpha,\alpha^2) \mid \alpha \in \fqssqq\}$,
$A=(\eta,b\eta^2)$. Then by Proposition \ref{laterale} the point
set $\cK':=\cK_{G}\cup \varphi_A(\cK_G)$ is a translation $(2\sqrt
q)$-arc. Note that half of the points of $\cK'$ are contained in
the conic of equation $X_2X_3=X_1^2$,  the other half in the conic
$X_2X_3=X_1^2+(b+1)\eta^2X_3^2$.
\end{example}

Example \ref{n3} provides hyperfocused arcs which are not
contained in any regular hyperoval. Actually,  the existence of
hyperfocused arcs which are neither contained in any hyperoval nor
in any proper subplane of $PG(2,q)$  can be proved. The proof of
such a result needs the following two lemmas.
\begin{lemma}\label{iper}
Let $\cK$ be a translation $q$-arc containing both points
$(0,0,1)$ and $(1,1,1)$. Then there exist $\alpha,\beta\in \fq$
and a positive integer $i$ with $(i,r)=1$, such that
$$\cK=\{(x,y,1)\mid \alpha x+(\alpha+1) y+\beta x^{2^i}+(\beta+1)
y^{2^i}=0\}\,.$$
\end{lemma}
\begin{proof}
Let $\psi_{\alpha\gamma}$ be the linear collineation
$$\psi_{\alpha,\gamma}:(X_1,X_2,X_3)\mapsto (\alpha X_1+(\alpha+1)X_2,\gamma
X_1+(\gamma+1)X_2,X_3)\,,$$ with $\alpha,\gamma\in \fq$, and let
$\cK'=\psi_{\alpha,\gamma}(\cK)$. Choose $\alpha$ and $\gamma$ in
such a way that the two points on $\ell_\infty$ which belong to no
secant of $\cK'$ are $(1,0,0)$ and $(0,1,0)$. Note that $\cK'$
contains $(0,0,1)$ and $(1,1,1)$. Also, for each $t\in \fq$ there
exists exactly one point $P_t$ of $\cK'$ on the line $X_2=t X_3$.
Let $F$ be the function on $\fq$ such that $P_t=(F(t),t,1)$. As
$\cK'$ is a translation arc containing $(0,0,1)$, the set
$\{(F(t),t)\mid t \in \fq\}$ is an additive subgroup of $\fq\times
\fq$. This implies $F(s+t)=F(s)+F(t)$ for any $s,t\in \fq$.
Theorem 8.41 in \cite{HIR} yields $F(t)=t^{2^i}$ for some $i$ with
$(i,r)=1$, that is
$$
\cK'=\{(x,y,1)\mid x=y^{2^i}\}\,,
$$
whence
$$\cK=\{(x,y,1)\mid (\alpha x+(\alpha+1)y)=(\gamma x+(\gamma+1)y)^{2^i}\}\,.$$
Then the assertion follows by letting $\beta=\gamma^{2^i}$.
\end{proof}

\begin{lemma}\label{noncop}
Assume that $r$ has a proper divisor $s>2$, and let $q'=2^s$. Let
$\cK=\cK_G$ with $G=\{(a,a^2)\mid a \in \fqss\}$. Then there exist
at most $\frac{r}{s}$ translation $q$-arcs containing $\cK$.
\end{lemma}
\begin{proof}
Let $\cI$ be any translation arc of size $q$ containing $\cK$.
Then by Lemma \ref{iper} there exist $\alpha,\beta\in \fq$ with
$\alpha\neq \beta$, and a positive integer $i$ with $(i,r)=1$,
such that
$$
\alpha a+(\alpha+1) a^2+\beta a^{2^i}+(\beta+1) a^{2^{i+1}}=0\,,
$$
for any $a\in \fqss$. This means that the polynomial $g(T):=
\alpha T+(\alpha+1) T^2+\beta T^{2^i}+(\beta+1) T^{2^{i+1}}$ must
be divisible by $T^{q'}+T$. If $2^{i+1}<q'$ this can only happen
for $g(T)\equiv 0$, that is $i=1$, $\beta=1$, $\alpha=0$. If
$2^{i+1}=q'$, that is $i=s-1$, then $\alpha=1$, $\beta=0$.
Finally, if $2^{i+1}>q'$, then also $2^i>q'$ as $(i,r)=1$. Write
$i=us+v$ with $u,v$ integers with $0\le v<s$. Then
$2^i=(q'-1)^{2^u}2^v+2^v$, and $g(T) \,\,{\rm mod}\,\, {T^{q'}+T}$
is the polynomial $\alpha T+(\alpha+1) T^2+\beta T^{2^v}+(\beta+1)
T^{2^{v+1}}$, which has to be the zero polynomial. This implies
$\beta=1$, $\alpha=0$,
$i\in\{s+1,2s+1,\ldots,(\frac{r}{s}-1)s+1\}$. Note that the arc
defined by $i=s-1$, $\alpha=1$, $\beta=0$ coincides with that
defined by $i=(\frac{r}{s}-1)s+1$, $\alpha=0$, $\beta=1$. Then the
assertion follows.
\end{proof}
Now we are in a position to prove the following theorem.
\begin{theorem}\label{sopra}
Let $q=2^r$ be such that there $r$ admits a proper divisor $s>2$.
Then  there exists a translation arc $\cK$ in $PG(2,q)$ such that
\begin{itemize}
\item[\rm{ (a)}] every point in $PG(2,q)\setminus \ell_\infty$
belongs to some secant of $\cK$;

\item[\rm{ (b)}] $\cK$ is not contained in any hyperoval;

\item[\rm{ (c)}] $\cK$ is not contained in any proper subplane.
\end{itemize}
\end{theorem}
\begin{proof}
Let $\cK_G$ be as in Lemma \ref{noncop}, and let $\cI_1,\ldots,
\cI_h$ be the translation $q$-arcs containing $\cK_G$. Note that
$h\le \frac{r}{s}$ by Lemma \ref{noncop}. As there are exactly
$q'(q'-1)/2$ secants of $\cK_G$, the number of  points in
$PG(2,q)\setminus \ell_\infty$ contained in no secant of $\cK_G$
is at least $q^2-q({q'}^2-q')/2= q(2^r-2^{2s-1}+2^{s-1})$. On the
other hand, the number of points in $\cup_{i=1,\ldots,h}\cI_i$ is
at most $\frac{qr}{s}$. It is straightforward to check that
$2^r-2^{2s-1}+2^{s-1}>\frac{r}{s}$. Hence, there exists a point
$\overline{A_1}\in PG(2,q)\setminus \ell_\infty$  which is
contained neither in a $\cI_i$ nor in a secant of $\cK_G$. Define
$G_1=G+A_1$ and $\cK_1=\cK_{G_1}$. If every point in $
PG(2,q)\setminus \ell_\infty$ belongs to some secant of $\cK_1$,
let $\cK:=\cK_1$. Otherwise choose a point $\overline{A_2}$ not
belonging to any secant of $\cK_1$ and let $G_2=G_1+A_2$,
$\cK_2=\cK_{G_2}$. Repeat the process until the arc $\cK_i$ has
the property that every point in $PG(2,q)\setminus \ell_\infty$
belongs to some secant of $\cK_i$, and define $\cK=\cK_i$.
Clearly, (a) is fulfilled by construction. Assume now that $\cK$
is contained in a hyperoval $\cI'$. By (a), $\cK$ coincides with
the points of $\cI'$ not on $\ell_\infty$, that is $\cK$ is one of
the translation $q$-arcs containing $\cK_G$. But this is
impossible as $\cK_{1}$ is not contained in any $\cI_i$ by
construction. Finally, (c) holds when $s$ is chosen to be the
maximum proper divisor of $r$. In fact, in this case the maximum
order of a subplane of $PG(2,q)$ is $2^s$, whereas $\#\cK\ge
2^{s+1}>2^s+2$.

\end{proof}

Theorem \ref{sopra} suggests that it might be hard to deal with
the problem of characterizing hyperfocused arcs.

\section{Generalized hyperfocused arcs}\label{s3}

In this section we consider generalized hyperfocused arcs, that is
arcs admitting a non-necessarily linear blocking set of minimum
size. In \cite{AKS} it is shown that an arc in $PG(2,q)$, $q$
even, does not admit a non-linear blocking set of its secants of
minimum size, provided that it is contained in a conic. The
following theorem proves that $k$-arcs admitting non-linear
blocking sets of size $k-1$ actually exist.

\begin{theorem}
Let $\cK$ be a translation $k$-arc, $k\ge4$, and let $\varphi$ be
a homology with axis $\ell_\infty$ and centre not in $\cK$. If the
set $\cK'=\cK\cup \varphi(\cK)$ is an arc, then there exists a
non-linear blocking set $\cB$ of the secants of $\cK'$ of minimum
size.
\end{theorem}

\begin{proof}
Assume that $(0,0,1)\in \cK$, and let $\cK=\cK_G$, with  $G$ an
additive subgroup of $\fq\times \fq$. Let $\overline{C}$ be the
centre of $\varphi$. Define $\cB$ as the subset of $2k-1$ points
$PG(2,q)$ which comprises  points $\overline{A}_\infty$, together
with the centres  of the homologies $\varphi \varphi_A$, with $A$
ranging over  $G\setminus \{(0,0)\}$. Let $l_{PQ}$ be any secant
of $\cK'$. If both $P$ and $Q$ are either in $\cK$ or in
$\varphi(\cK)$, then $l_{PQ}$ meets $\cB$ in a point
$\overline{A}_\infty$, for some $A \in G\setminus \{(0,0)\}$. Now
assume that $P=\overline{A}$ and $Q=\varphi(\overline{B})$ for
some $A,B\in G$. Then $l_{PQ}$ passes through the centre of
$\varphi\varphi_{A+B}$. This proves that $\cB$ is a blocking set
of the secants of $\cK'$. As $\cB$ has size $2k-1$ and is not
contained in any line, the assertion is proved.
\end{proof}

\begin{example}\label{otto}
Let $\cK=\cK_G$ with $G=\{(0,0),(0,1),(1,0),(1,1)\}$. Consider the
homology

\begin{equation}\label{otto1}
 \varphi:(X_1,X_2,X_3)\mapsto (\lambda X_1+a_1X_3,\lambda X_2
+a_2X_3,X_3)\,,
\end{equation}
 with
\begin{itemize}
\item $\lambda \in \fq$, $\lambda\neq 0,1$, $a_1,a_2\in \fq$;

\item $\{a_1,a_2,a_1+a_2\}\cap
\{0,1,\lambda,\lambda+1\}=\emptyset$.

\end{itemize}

Then it is straightforward to check that $\cK'=\cK\cup
\varphi(\cK)$ is an arc. A non-linear blocking set $\cB$ of the
secants of $\cK'$ of minimum size is
$$
\begin{array}{cl}
\cB=& \{(1,0,0),(0,1,0),(1,1,0),(a_1,a_2,1+\lambda),
(a_1+\lambda,a_2,1+\lambda),\\{}&
(a_1,a_2+\lambda,1+\lambda),(a_1+\lambda,a_2+\lambda,1+\lambda)\}\,,
\end{array}
$$
which consists of the points of a subplane of $PG(2,q)$ of order
$2$.
\end{example}

The following result shows that a non-linear blocking set of
minimum size of the secants of a $k$-arc cannot be an arc itself.
Also, it will be useful for the classification of small
generalized hyperfocused arcs which will be given in next section.

\begin{proposition}\label{triangle}
Let $\cB$ be a blocking set of minimum size of the secants of a
$k$-arc $\cK$ in $PG(2,q)$, $q$ even. Then any three points in
$\cB$ blocking the secants of a $3$-arc contained in $\cK$ are
collinear.
\end{proposition}
\begin{proof}
This proof relies on the idea of  Segre's celebrated Lemma of
Tangents \cite{SEG}. Let  $P_1$, $P_2$ and $P_3$ be any three
distinct points in $\cK$. For each $i\in\{1,2,3\}$, let $Q_i\in
\cB$ be collinear with $P_{j}$ and $P_{k}$, where $j,k\in
\{1,2,3\}$, $j=i+1 \pmod 3$, $k=i-1 \pmod 3$. It has to be proved
that $Q_1$, $Q_2$ and $Q_3$ are collinear. Assume without loss of
generality that $P_1=(1,0,0)$, $P_2=(0,1,0)$ and $P_3=(0,0,1)$.
For a point $P$ distinct from $P_i$, $i=1,2,3$, let $\alpha_P^1$,
$\alpha_P^2$, $\alpha_P^3$ be the elements of $\fq$ such that
\begin{itemize}
\item $X_3=\alpha_P^1 X_2$ is the line through $P_1$ and $P$,

\item $X_1=\alpha_P^2 X_3$ is the line through $P_2$ and $P$,

\item $X_2=\alpha_P^3 X_1$ is the line through $P_3$ and $P$.
\end{itemize}
It is straightforward to check that if $P$ does not belong to the
triangle with vertices $P_1$, $P_2$, $P_3$, then
\begin{equation}\label{uno}
\alpha_P^1\alpha_P^2\alpha_P^3=1\,.
\end{equation}
Now, consider the set of secants of $\cK$ passing through exactly
one point among $P_1$, $P_2$ and $P_3$. Clearly, it coincides with
the set which comprises the lines joining $P_1$, $P_2$ and $P_3$
to any point of $\cB\setminus \{Q_1,Q_2,Q_3\}$, together with the
lines through $P_i$ and $Q_i$, $i=1,2,3$. Hence,
$$
\prod_{P\in \cK,\,\,P\neq P_1,P_2,P_3}
\alpha_P^1\alpha_P^2\alpha_P^3=\alpha_{Q_1}^1\alpha_{Q_2}^2\alpha_{Q_3}^3\left(\prod_{Q\in
\cB,\,\,Q\neq Q_1,Q_2,Q_3}
\alpha_Q^1\alpha_Q^2\alpha_Q^3\right)\,.
$$
Then by (\ref{uno}),
$\alpha_{Q_1}^1\alpha_{Q_2}^2\alpha_{Q_3}^3=1$ holds. As $q$ is
even, this is equivalent to the collinearity of $Q_1$, $Q_2$ and
$Q_3$ and the assertion is proved.
\end{proof}

\section{Classification of small generalized hyperfocused arcs}\label{s4}
The aim of this section if to classify the small arcs admitting
blocking sets of minimum size for their secants. The linear case
has already been settled in \cite{DK} and \cite{CH}. The main
result of the section is the following.
\begin{theorem}\label{classify1}
Let $\cK$ be a $k$-arc in $PG(2,q)$, $q$ even, with $k\le 10$. If
there exists a minimal non-linear blocking set of the secants of
$\cK$, then $k=8$ and $\cK$ is projectively equivalent to the arc
$\cK'$ in Example \ref{otto}.
\end{theorem}

The proof of this result relies on a connection between blocking
sets of the secants of an arc and $1$-factorizations of complete
graphs. For the sake of completeness, some basic definitions from
graph theory are reported.

Let $K_{2n}$ be the complete graph with $2n$ vertices. A {\em
$1$-factor} of $K_{2n}$ is a set of vertex disjoint edges which
cover the vertices of $K_{2n}$. An edge disjoint set of
$1$-factors covering the edges of $K_{2n}$ is said to be a {\em
$1$-factorization} of $K_{2n}$. The set of vertices of $K_{2n}$
will be denoted by $V(K_{2n})$.

\begin{definition}\label{embe} Let $\cF$ be a $1$-factorization of $K_{2n}$.
An {\em embedding} of $\cF$ in $PG(2,q)$ is an injective map
$\psi:V(K_{2n})\cup \cF\rightarrow PG(2,q)$ such that
\begin{enumerate}
\item[\rm{i) }] for any $i,j,k \in V(K_{2n})$, the points
$\psi(i)$, $\psi(j)$, $\psi(k)$ are not collinear;

\item[\rm{ii) }] for any $F \in \cF$, the point $\psi(F)$ is
collinear with $\psi(i)$ and $\psi(j)$, for every edge $(i,j) \in
F$.
\end{enumerate}
\end{definition}
Given an embedding $\psi$ of a $1$-factorization $\cF$ of $K_{2n}$
in $PG(2,q)$, the set $\psi(V(K_{2n}))$ is an arc, whereas
$\psi(\cF)$ is a blocking set of minimum size of the secant of
such arcs. The following equivalent formulation of Theorem
\ref{classify1} will be proved.
\begin{theorem}\label{classify}
Let $\psi$ be an embedding of a $1$-factorization $\cF$ of
$K_{2n}$ in $PG(2,q)$, $q$ even, with $3\le n\le 5$. If the points
$\{\psi(F)\mid F \in \cF\}$ are not collinear, then $n=4$ and
$\psi(V(K_{2n}))$ is projectively equivalent to the arc $\cK'$  in
Example \ref{otto}.
\end{theorem}

Assume that $V(K_{2n})=\{1,2,\ldots,2n\}$, $n \geq 3$, and let
$\cF=\{F_1,F_2,\ldots,F_{2n-1}\}$ be a $1$-factorization of
$K_{2n}$. Let $\psi$ be an embedding of $\cF$ in $PG(2,q)$.

\subsection{Proof of Theorem \ref{classify} for $n=3$}
As all the $1$-factorizations of the complete graph with $6$
vertices are isomorphic,  we may assume that:
\begin{itemize}
\item $\psi(F_1)$ is the common point of the lines
$\psi(1)\psi(2)$, $\psi(3)\psi(4)$, $\psi(5)\psi(6)$;

\item $\psi(F_2)$ is the common point of the lines
$\psi(1)\psi(3)$, $\psi(2)\psi(5)$, $\psi(4)\psi(6)$;

\item $\psi(F_3)$ is the common point of the lines
$\psi(1)\psi(4)$, $\psi(2)\psi(6)$, $\psi(3)\psi(5)$;

\item$\psi(F_4)$ is the common point of the lines
$\psi(1)\psi(5)$, $\psi(2)\psi(4)$, $\psi(3)\psi(6)$;

\item $\psi(F_5)$ is the common point of the lines
$\psi(1)\psi(6)$, $\psi(2)\psi(3)$, $\psi(4)\psi(5)$.
\end{itemize}

By Proposition \ref{triangle} the following triples of points are
collinear:
$$ \psi(F_1),\psi(F_2),\psi(F_3), \qquad
\psi(F_1),\psi(F_2),\psi(F_4),\qquad
\psi(F_1),\psi(F_2),\psi(F_5)\,.$$ Then  all points in
$\{\psi(F)\mid F\in \cF\}$ are collinear, which proves the
assertion.

\subsection{Proof of Theorem \ref{classify} for $n=4$}
There are $6$ non-isomorphic $1$-factorizations of $K_8$ (see e.g.
\cite{LARS}). From the proof of Theorem 5.3 in \cite{CH},
it follows that $4$ of them
cannot be embedded in $PG(2,q)$. We are  left with the following
two cases.

{\bf Case 1}: $\cF=\{F_1,\ldots,F_7\}$ with
$$
\begin{array}{ll}
F_1=\{(8,1),(2,3),(4,5),(6,7)\},&
F_2=\{(8,2),(1,3),(4,6),(5,7)\},\\
F_3=\{(8,3),(1,2),(4,7),(5,6)\},&
F_4=\{(8,4),(1,5),(2,6),(3,7)\},\\
F_5=\{(8,5),(1,4),(2,7),(3,6)\},&
F_6=\{(8,6),(1,7),(2,4),(3,5)\},\\
F_7=\{(8,7),(1,6),(2,5),(3,4)\}.&{}
\end{array}$$
Assume without loss of generality that $\psi(4)=(0,0,1)$,
$\psi(5)=(0,1,1)$, $\psi(6)=(1,0,1)$, $\psi(7)=(1,1,1)$, that is
$\{\psi(4),\psi(5),\psi(6),\psi(7)\}$ coincides with $\cK_G$, with
$G$ as in Example \ref{otto}. Then $\psi(F_1)=(0,1,0)$,
$\psi(F_2)=(1,0,0)$ and $\psi(F_3)=(1,1,0)$. Now, note that by
Proposition \ref{triangle} the following triples of points are
collinear: $$ \psi(F_4),\psi(F_5),\psi(F_1), \qquad
\psi(F_4),\psi(F_6),\psi(F_2),\qquad
\psi(F_4),\psi(F_7),\psi(F_3)\,.$$ Hence, if $\psi(F_4)$ lies on
$\ell_\infty$, then the whole $\{\psi(F)\mid F \in \cF\}$ is
contained in a line. Now assume that $\psi(F_4)\notin
\ell_\infty$.

Let $\varphi$ be the linear collineation of $PG(2,q)$ such that
$\varphi(\psi(4))=\psi(8)$, $\varphi(\psi(5))=\psi(1)$,
$\varphi(\psi(6))=\psi(2)$ and $\varphi(\psi(7))=\psi(3)$.
Clearly, $\varphi$ fixes $\psi(F_1)$, $\psi(F_2)$, $\psi(F_3)$,
and hence $\varphi$ is a central collineation with axis
$\ell_\infty$. The centre of $\varphi$ is  $\psi(F_4)$, which is
assumed not to belong to $\ell_\infty$. Therefore, $\varphi$ is as
in Equation ({\ref{otto1}) for some $a_1,a_2,\lambda\in \fq$. As
$\psi(V(K_8))=\cK_G\cup \varphi(\cK_G)$ is an arc, it is
straightforward to check that
\begin{itemize}
\item $\lambda \in \fq$, $\lambda\neq 0,1$, $a_1,a_2\in \fq$;
\item $\{a_1,a_2,a_1+a_2\}\cap
\{0,1,\lambda,\lambda+1\}=\emptyset$.
\end{itemize}
Then the assertion is proved.

 {\bf Case 2}:
 $\cF=\{F_1,\ldots,F_7\}$ with
$$
\begin{array}{ll}
F_1=\{(8,1),(2,3),(4,5),(6,7)\},&
F_2=\{(8,2),(1,4),(3,6),(5,7)\},\\
F_3=\{(8,3),(1,6),(2,5),(4,7)\},&
F_4=\{(8,4),(1,7),(2,6),(3,5)\},\\
F_5=\{(8,5),(1,2),(3,7),(4,6)\},&
F_6=\{(8,6),(1,5),(2,7),(3,4)\},\\
F_7=\{(8,7),(1,3),(2,4),(5,6)\}.&{}
\end{array}
$$
By Proposition \ref{triangle},
 any point $\psi(F_i)$ with $3\le
i\le 7 $ is collinear with $\psi(F_1)$ and $\psi(F_2)$. Then all
points in $\{\psi(F)\mid F\in \cF\}$ are collinear, which proves
the assertion.

\subsection{Proof of Theorem \ref{classify} for $n=5$}
Define $\cT_\cF^{0}$ as the set of all triples $\{F_i,F_j,F_k\}$
such that $(i,j)\in F_k$, $(i,k)\in F_j$, $(j,k)\in F_i$, with
$i,j,k$ ranging over $V(K_{10})$. By Proposition \ref{triangle},
for any $\{F_i,F_j,F_k\}\in \cT_\cF^{0}$ the points $\psi(F_i)$,
$\psi(F_j)$ and $\psi(F_k)$ are collinear.

Now define recursively a set $\cT_\cF^i$, $i\ge 1$, as follows:
$\cT_\cF^i$ contains all the joins of two sets in $\cT_\cF^{i-1}$
sharing at least two elements of $\cF$. Clearly, for any $\cA\in
\cT_\cF^{i}$, the points $\{\psi(F)\mid F \in \cA\}$ are
collinear. By the following lemma, all points in $\{\psi(F)\mid F
\in \cF\}$ are collinear, which completes the proof of Theorem
\ref{classify}.

\begin{lemma}\label{computer}
For any $1$-factorization $\cF$ of $K_{10}$, there exists an
integer $i$ for which $\cT_\cF^{i}$ contains $\cF$.
\end{lemma}
The proof of Lemma \ref{computer} consists of a computer based
investigation of all 396 non-isomorphic $1$-factorizations of
$K_{10}$ (\cite{LARS}). For the details of the proof the reader is
referred to \cite{GM}.

\end{document}